\documentclass[12pt]{amsart}

\usepackage{amsmath,amsfonts,amscd,amsthm,amssymb,latexsym}

\font\teneufm=eufm10 \font\seveneufm=eufm7 \font\fiveeufm=eufm5
\newfam\frakturfam
\textfont\frakturfam=\teneufm \scriptfont\frakturfam=\seveneufm
\scriptscriptfont\frakturfam=\fiveeufm

\newtheorem{lm}{Lemma}
\newtheorem{theor}{Theorem}
\newtheorem{co}{Corollary}

\def\bee{\begin{eqnarray}}
\def\bes{\begin{eqnarray*}}
\def\eee{\end{eqnarray}}
\def\ees{\end{eqnarray*}}

\def\a{\alpha}
\def\b{\beta}

\def\Proof{{\sl Proof.}\ }

\title{Automorphisms of Veronese subalgebras of polynomial algebras and free Poisson algebras}

\begin{document}
\date{}
\maketitle

\begin{center}
{\bf Bakhyt Aitzhanova}\footnote{Department of Mathematics,
 Wayne State University,
Detroit, MI 48202, USA, e-mail: {\em aitzhanova.bakhyt01@gmail.com}}, 
{\bf Leonid Makar-Limanov}\footnote{Department of Mathematics, Wayne State University, Detroit, MI 48202, USA 
and Faculty of Mathematics and Computer Science, Weizmann Institute of Science, Rehovot, 7610001, Israel, 
e-mail: {\em lml@wayne.edu}}, 
and
{\bf Ualbai Umirbaev}\footnote{Department of Mathematics,
 Wayne State University,
Detroit, MI 48202, USA 
and Institute of Mathematics and Mathematical Modeling, Almaty, 050010, Kazakhstan,
e-mail: {\em umirbaev@wayne.edu}}

\end{center}

\begin{abstract} The Veronese subalgebra $A_0$ of degree $d\geq 2$ of the polynomial algebra $A=K[x_1,x_2,\ldots,x_n]$ over a field $K$ in the variables $x_1,x_2,\ldots,x_n$ is the subalgebra of $A$ generated by all monomials of degree $d$ and the Veronese subalgebra $P_0$ of degree $d\geq 2$ of the free Poisson algebra $P=P\langle x_1,x_2,\ldots,x_n\rangle$ is the subalgebra spanned by all homogeneous elements of degree $kd$, where $k\geq 0$. 

If $n\geq 2$ then every derivation and every locally nilpotent derivation of $A_0$ and $P_0$ over a field $K$ of characteristic zero is induced by a derivation and a locally nilpotent derivation of $A$ and $P$, respectively. Moreover, we prove that every automorphism of $A_0$ and $P_0$ over a field $K$ closed with respect to taking all $d$-roots of elements is induced by an automorphism of $A$ and $P$, respectively. 
\end{abstract}


\noindent


\section{Introduction}

\hspace*{\parindent}

Let $K$ be an arbitrary field and let $\mathbb{A}^n$ and $\mathbb{P}^n$ be the affine and the projective $n$-spaces over $K$, respectively. 
The \textit{Veronese map} of degree $d$ is the map
\bes
\nu_{n,d}:\mathbb{P}^n\to\mathbb{P}^m
\ees 
that sends $[x_0:\ldots:x_n]$ to all $m+1$ possible monomials of total degree $d$, where 
\bes
m={n+d \choose d}-1.
\ees
It is well known that the image $V_{n,d}$ of the Veronese map $\nu_{n,d}$ is a projective variety and is called the \textit{Veronese variety} \cite{Har}. 

If $Y_{i_0\ldots i_n}=x_0^{i_0}\ldots x_n^{i_n}$, $i_0+\ldots+i_n=d$, then the Veronese variety is determined by the set of quadratic relations 
\bee\label{h1}
Y_{i_0\ldots i_n}Y_{j_0\ldots j_n}=Y_{k_0\ldots k_n}Y_{r_0\ldots r_n},
\eee
where $i_0+j_0=k_0+r_0,\ldots,i_n+j_n=k_n+r_n$ \cite{Shaf}.

The {\em affine cone of $V_{n,d}$} or the {\em affine Veronese variety} is the affine subvariety $\mathbb{A}^{m+1}$ determined by the set of relations (\ref{h1}). The algebra of polynomial functions on the affine Veronese cone  $V_{n,d}$ is isomorphic to the subalgebra of $K[x_0,\ldots,x_n]$ generated by all monomials of degree $d$ \cite{TGI}. 

  Veronese surfaces play an important role in the description of quasihomogeneous affine surfaces given by M.H. Gizatullin \cite{Giz1971} and V.L. Popov \cite{Popov1973}. They form one of the main examples of the so-called {\em Gizatullin surfaces} \cite{KPZ2017}. The structure of the automorphism groups of Veronese surfaces are studied in \cite{AU,AZ2013,GD1975,GD1977,ML90,ML01}. The derivations and locally nilpotent derivations of affine Veronese surfaces are described in \cite{AU}.

Recently J. Kollar devoted two interesting papers \cite{Kollar1, Kollar2} to the study of automorphism groups of some more general affine varieties. In particular, he described the group of automorphisms of all affine Veronese varieties.

Let $A=K[x_1,x_2,\ldots,x_n]$ be the polynomial algebra over a field $K$ in the variables $x_1,x_2,\ldots,x_n$. Consider the 
grading 
\bes
A=A_0\oplus A_1\oplus \ldots \oplus A_{d-1}, 
\ees
where $d\geq 2$ and $A_i$ is the subspace of $A$ generated by all monomials of degree $kd+i$ for all $k\geq 0$.  This is a $\mathbb{Z}_d$-grading of $A$, i.e., $A_iA_j\subseteq A_{i+j}$ for all $i,j\in \mathbb{Z}_d$. The subalgebra $A_0$ is called the {\em Veronese subalgera of $A$ of degree $d$}. 

Similarly, let $P=P\langle x_1,x_2,\ldots,x_n\rangle$ be the free Poisson algebra over $K$ in the variables $x_1,x_2,\ldots,x_n$. The grading 
\bes
P=P_0\oplus P_1\oplus \ldots \oplus P_{d-1}, 
\ees
where $P_i$ is the subspace of $P$ generated by all homogeneous elements of degree $kd+i$ is a $\mathbb{Z}_d$-grading of $P$, i.e., $P_iP_j\subseteq P_{i+j}$ and $\{P_i,P_j\}\subseteq P_{i+j}$ for all $i,j\in \mathbb{Z}_d$. The subalgebra $P_0$ is called the {\em Veronese subalgera of $P$ of degree $d$}. 

We prove that every derivation and every locally nilpotent derivation of $P_0$ (of $A_0$) over a field $K$ of characteristic zero is induced by a derivation and a locally nilpotent derivation of $P$ (of $A$), respectively. Moreover, we prove that every automorphism of $P_0$ (of $A_0$) over a field $K$, closed with respect to taking all $d$-roots of elements, is induced by an automorphism of $P$ (of $A$). 

This paper is organized as follows. In Section 2 we describe a basis of free Poisson algebras and give some elementary definitions that are necessary in the future. Derivations of the Veronese subalgebras are studied in Section 3 and automorphisms are studied in Section 4. All results are proven for Poisson algebras. Analogues of these results for polynomial algebras are formulated in Section 5. In the same section we give some counter-examples that show the analogous results do not hold for polynomial algebras in one variable and for free associative algebras.

\section{Free Poisson algebra $P\langle x_1,\ldots,x_n\rangle$}

\hspace*{\parindent}

A vector space $P$ over a field $K$ endowed with two bilinear operations $x\cdot y$  (a multiplication) and $\{x,y\}$ (a Poisson bracket) is called a \emph{Poisson algebra} if $P$ is a commutative associative algebra under $x\cdot y$, $P$ is a Lie algebra under  $\{x,y\}$, and $P$ satisfies the following identity (the Leibniz identity):
$$ \{x,y\cdot z\}=y\cdot \{x,z\}+ \{x,y\}\cdot z.$$

Symplectic Poisson algebras $P_n$ appear in many areas of algebra. For each natural $n\geq 1$ the symplectic Poisson algebra $P_n$ of index $n$ is a polynomial algebra
$K[x_1,y_1, \ldots,x_n,y_n]$ endowed with the Poisson bracket
defined by
\bes \{x_i,y_j\}=\delta_{ij}, \ \ \{x_i,x_j\}=0, \ \
\{y_i,y_j\}=0, \ees
where $\delta_{ij}$ is the Kronecker symbol
and $1\leq i,j\leq n$.

Let $P$ be a Poisson algebra. A linear map $D: P\to P$ is called a {\em derivation} of the Poisson algebra $P$ if
\bee\label{h2}
D(xy)=D(x)y+xD(y)
\eee
and
\bee\label{h3}
D(\{x,y\})=\{D(x),y\}+\{x,D(y)\}
\eee
for all $x,y\in P$. 

We refer to a linear map $D: P\to P$ satisfying (\ref{h2}) as an {\em associative derivation} of $P$ and to a linear map $D: P\to P$ satisfying (\ref{h3}) as a {\em Lie derivation} of $P$. Thus, a derivation of a Poisson algebra is simultaneously an associative and a Lie derivation. We often refer to them as Poisson derivations. 

The Leibniz  identity implies that for any $x\in P$ the map
\bes 
\mathrm{ad}_x: P\to P \,\,\, (y\mapsto \{x,y\})
\ees
is an associative derivation of $P$. 
The Jacobi identity implies that the map $\mathrm{ad}_x$ is also a Lie derivation of $P$. 
 So for any $x\in P$ the map $\mathrm{ad}_x$ is a Poisson derivation of $P$.

Let  $L$  be a Lie algebra with Lie bracket $[\,,]$ over a field $K$ and let
$e_1,e_2\ldots,e_k,\ldots$ be a linear basis of $L$. Then there exists a unique bracket $\{\,,\}$ on the
 polynomial algebra $K[e_1, e_2,\ldots,e_k,\ldots]$  defined by 
\bes
\{e_i,e_j\}=[e_i,e_j]
\ees
 for all $i,j$ and satisfying the Leibniz identity.  
With this bracket   
\bes
P(L)=\langle K[e_1,e_2,\ldots],\cdot, \{\,,\}\rangle
\ees
 becomes a Poisson algebra. This Poisson algebra $P(L)$ is called the {\em Poisson enveloping algebra} \cite{UZh20} (or {\em Poisson symmetric algebra} \cite{MLU11JA}) of $L$. Note that the bracket $\{,\}$ of the algebra $P(L)$ depends on the structure of $L$ but does not depend on a chosen basis.  

Let $L=\mathrm{Lie}\langle x_1,\ldots, x_n \rangle$ be the free Lie algebra with free generators $x_1,\ldots, x_n$. It is well-known (see, for example \cite{Sh}) that $P(L)$ is the free Poisson algebra over $K$ in the variables $x_1,\ldots, x_n$. We denote it by $P\langle x_1,\ldots, x_n\rangle$.  

Let us choose a multihomogeneous linear basis
\bes 
x_1, \ldots, x_n, [x_1,x_2],\ldots,[x_1,x_n],\ldots, [x_{n-1},x_n], [[x_1,x_2],x_3], \ldots
\ees
of a free Lie algebra $L$ and denote the elements of this basis by 
\bee\label{h4}
e_1,e_2,\ldots,e_s,\ldots.
\eee
The algebra $P=P\langle x_1, \ldots, x_n\rangle$ coincides with the polynomial algebra on the elements (\ref{h4}). Consequently, the monomials 
\bee\label{h5}
u=e_{i_1}e_{i_2}\ldots e_{i_s},
\eee
where $i_1\leq i_2\leq\ldots\leq i_s$ form a linear basis of $P$.

Denote by $\deg$ the standard degree function on $P$, i.e., $\deg(x_i)=1$ for all $1\leq i\leq n$. 
If $u$ is an element of the form (\ref{h5}) then 
\bes
\deg u=\deg e_{i_1}+\deg e_{i_2}+\ldots +\deg e_{i_s}. 
\ees  
Set also $d(u)=s$ and call it the {\em polynomial length} 
of $u$. 
 Note that
\bes
\deg \{f,g\}=\deg f+\deg g
\ees
if $f$ and $g$ are homogeneous and $\{f,g\}\neq0$. 

Denote by $Q(P)=P(x_1,\ldots,x_n)$ the field of fractions of the polynomial algebra $K[e_1,e_2,\ldots]$ in the variables $(\ref{h4})$. The Poisson bracket $\{\cdot,\cdot\}$ on $K[e_1,e_2,\ldots]=P$ can be uniquely extended to a Poisson bracket on the field of its fractions $Q(P)$ 
and
\bee\label{h6}
\Bigl\{\frac{a}{b},\frac{c}{d}\Bigr\}=\frac{\{a,c\}bd-\{a,d\}bc-\{b,c\}ad+\{b,d\}ac}{b^2d^2}
\eee
for all $a,b,c,d\in P$ with $bd\neq 0$.

The field $Q(P)=P(x_1,x_2,\ldots,x_n)$ with this Poisson bracket is called the {\em free Poisson field} over $K$ in variables $x_1,\ldots, x_n$ \cite{MLSh}.

Several combinatorial results on the structure of free Poisson algebras and free Poisson fields are proven in \cite{MLU07,MLTU09JA,MLU11JA,MLSh,MLU16,12JA}.

We fix a grading
\bee\label{h7}
P=P_0\oplus P_1\oplus \ldots\oplus P_{d-1}
\eee
of the free Poisson algebra $P=P\langle x_1,\ldots,x_n\rangle$, where $P_i$ is the linear span of all elements of degree $i+ds$, $i=0,1,\ldots,d-1$, and $s$ is an arbitrary nonnegative integer. This is a $\mathbb{Z}_d$-grading of $P$, i.e., 
\bes
P_iP_j\subseteq P_{i+j}, \ \ \{P_i,P_j\} \subseteq P_{i+j},
\ees
 where $i,j\in  \mathbb{Z}_d=\mathbb{Z}/d\mathbb{Z}$. For shortness we will refer to this grading as the $d$-grading.

An automorphism $\phi\in\mathrm{Aut}\,P$ is called a {\em graded automorphism} with respect to grading  (\ref{h7}) if $\phi(x_1), \phi(x_2),\ldots,\phi(x_n)\in P_1$. A graded automorphism is called {\em graded tame} if it is a product of graded elementary automorphisms. 

We will call a graded automorphism  of $P$ with respect to grading  (\ref{h7}) a {\em d-graded automorphism} for shortness.  Obviously, every $d$-graded automorphism induces an automorphism of the algebra $P_0$. 
A  derivation $D$ of $P$ will be called a {\em d-graded derivation} if 
 $D(x_1), D(x_2),\ldots,D(x_n)\in P_1$.

\section{Derivations of $P_0$}

\hspace*{\parindent}

In this section we assume that $K$ is an arbitrary field of characteristic zero. 
\begin{lm}\label{m1} Every derivation of the Poisson algebra $P=P\langle x_1,x_2,\ldots,x_n\rangle$ over $K$ can be uniquely extended to a derivation of the Poisson field $Q(P)=P(x_1,x_2,\ldots,x_n)$. 
\end{lm}
\Proof Let $D$ be an arbitrary derivation of the free Poisson algebra $P$. In particular, $D$ 
is a derivation of the polynomial algebra $P=K[e_1,\ldots,e_s,\ldots]$. It is well known \cite[p. 120]{ZS} that $D$ can be uniquely extended to an associative  derivation $S$ of the quotient field $Q(P)=K(e_1,\ldots,e_s,\ldots)$. 
We will show that $S$ is a Lie derivation of $Q(P)$, i.e., that $S$ satisfies (\ref{h3}). We have to check that 
\bee\label{h8}
S\left(\Bigl\{\frac{a}{b},\frac{c}{d}\Bigr\}\right)=\Bigl\{S\left(\frac{a}{b}\right),\frac{c}{d}\Bigr\}+\Bigl\{\frac{a}{b},S\left(\frac{c}{d}\right)\Bigr\}
\eee
for all $a,b,c,d\in P\langle x_1,\dots,x_n\rangle$ with $bd\neq0$. 
Using (\ref{h6}) we get 
\bes
S\left(\Bigl\{\frac{a}{b},\frac{c}{d}\Bigr\}\right)=S\left(\frac{\{a,c\}bd-\{a,d\}bc-\{b,c\}ad+\{b,d\}ac}{b^2d^2}\right)\\
=S\left(\frac{\{a,c\}}{bd}\right)
-S\left(\frac{\{a,d\}c}{bd^2}\right)
-S\left(\frac{\{b,c\}a}{b^2d}\right)+S\left(\frac{\{b,d\}ac}{b^2d^2}\right)\\
=\frac{S(\{a,c\})bd-\{a,c\}S(bd)}{b^2d^2}
-\frac{S(\{a,d\}c)bd^2-\{a,d\}cS(bd^2)}{b^2d^4}\\
-\frac{S(\{b,c\}a)b^2d-\{b,c\}aS(b^2d)}{b^4d^2}
+\frac{S(\{b,d\}ac)b^2d^2-\{b,d\}acS(b^2d^2)}{b^4d^4}\\
=\frac{S(\{a,c\})}{bd}-\frac{\{a,c\}S(b)}{b^2d}-\frac{\{a,c\}S(d)}{bd^2}
-\frac{S(\{a,d\})c}{bd^2}-\frac{\{a,d\}S(c)}{bd^2}+\frac{\{a,d\}cS(b)}{b^2d^2}\\+\frac{2\{a,d\}cS(d)}{bd^3}-\frac{S(\{b,c\})a}{b^2d}-\frac{\{b,c\}S(a)}{b^2d}
+\frac{2\{b,c\}aS(b)}{b^3d}+\frac{\{b,c\}aS(d)}{b^2d^2}\\
+\frac{S(\{b,d\})ac}{b^2d^2}+\frac{\{b,d\}S(a)c}{b^2d^2}+\frac{\{b,d\}aS(c)}{b^2d^2}-\frac{2\{b,d\}acS(b)}{b^3d^2}
-\frac{2\{b,d\}acS(d)}{b^2d^3}\\
=\frac{\{S(a),c\}}{bd}+\frac{\{a,S(c)\}}{bd}
-\frac{\{a,c\}S(b)}{b^2d}-\frac{\{a,c\}S(d)}{bd^2}-\frac{\{S(a),d\}c}{bd^2}\\
-\frac{\{a,S(d)\}c}{bd^2}-\frac{\{a,d\}S(c)}{bd^2}+\frac{\{a,d\}cS(b)}{b^2d^2}+\frac{2\{a,d\}cS(d)}{bd^3}-\frac{\{S(b),c\}a}{b^2d}\\-\frac{\{b,S(c)\}a}{b^2d}
-\frac{\{b,c\}S(a)}{b^2d}+\frac{2\{b,c\}aS(b)}{b^3d}
+\frac{\{b,c\}aS(d)}{b^2d^2}+\frac{\{S(b),d\}ac}{b^2d^2}\\+\frac{\{b,S(d)\}ac}{b^2d^2}
+\frac{\{b,d\}S(a)c}{b^2d^2}
+\frac{\{b,d\}aS(c)}{b^2d^2}
-\frac{2\{b,d\}acS(b)}{b^3d^2}-\frac{2\{b,d\}acS(d)}{b^2d^3}.
\ees
Direct calculations give that 
\bes
\Bigl\{S\left(\frac{a}{b}\right),\frac{c}{d}\Bigr\}+\Bigl\{\frac{a}{b},S\left(\frac{c}{d}\right)\Bigr\}=\Bigl\{\frac{S(a)b-aS(b)}{b^2},\frac{c}{d}\Bigr\}+\Bigl\{\frac{a}{b},\frac{S(c)d-cS(d)}{d^2}\Bigr\}\\
=\Bigl\{\frac{S(a)}{b},\frac{c}{d}\Bigr\}-\Bigl\{\frac{aS(b)}{b^2},\frac{c}{d}\Bigr\}+\Bigl\{\frac{a}{b},\frac{S(c)}{d}\Bigr\}-\Bigl\{\frac{a}{b},\frac{cS(d)}{d^2}\Bigr\}\\
=\frac{\{S(a),c\}}{bd}
-\frac{\{S(a),d\}c}{bd^2}-\frac{\{b,c\}S(a)}{b^2d}+\frac{\{b,d\}S(a)c}{b^2d^2}-\frac{\{a,c\}S(b)}{b^2d}\\-\frac{\{S(b),c\}a}{b^2d}
+\frac{\{a,d\}S(b)c}{b^2d^2}
+\frac{\{S(b),d\}ac}{b^2d^2}+\frac{2\{b,c\}aS(b)}{b^3d}-\frac{2\{b,d\}aS(b)c}{b^3d^2}\\+\frac{\{a,S(c)\}}{bd}-\frac{\{a,d\}S(c)}{bd^2}-\frac{\{b,S(c)\}a}{b^2d}+\frac{\{b,d\}aS(c)}{b^2d^2}-\frac{\{a,c\}S(d)}{bd^2}\\-\frac{\{a,S(d)\}c}{bd^2}+\frac{2\{a,d\}cS(d)}{bd^3}+\frac{\{b,c\}S(d)a}{b^2d^2}
+\frac{\{b,S(d)\}ac}{b^2d^2}-\frac{2\{b,d\}acS(d)}{b^2d^3}.
\ees

These two equalities imply (\ref{h8}). $\Box$

Consider the grading (\ref{h7})  of $P$. 
A  Poisson derivation $D$ of $P$ will be called a {\em d-graded Poisson derivation} if $D(x_i)\in P_1$ for all $i=1,\ldots,n$. 
Obviously, every $d$-graded Poisson derivation of $P$ induces a Poisson derivation of 
$P_0$. The reverse is also true.

\begin{lm}\label{m2}
Every Poisson derivation of $P_0$ can be uniquely extended to a $d$-graded Poisson derivation of $P=P\langle x_1,x_2,\ldots,x_n\rangle$. 
\end{lm}
\Proof Let $D$ be a Poisson derivation of $P_0$. In particular, $D$ is a derivation of the associative and commutative algebra 
$P_0$. Since $P_0$ is a domain, $D$ can be uniquely extended \cite[p. 120]{ZS}
 to a derivation $T$ of the field of fractions $Q(P_0)$ of $P_0$. The field extension 
$$Q(P_0)\subseteq Q(P)$$
 is algebraic since every $e_i$ is a root of the polynomial $p(t)=t^d-e_i^d\in Q(P_0)[t]$ for all $i$. This extension is separable since $K$ is a field of characteristic zero. By Corollaries $2$ and $2'$ in \cite[pages 124--125]{ZS}, the associative derivation $T$ of the field 
$Q(P_0)$ can be uniquely extended to an associative derivation $S$ of the field $Q(P)$. 

Suppose that 
\bes
S(e_j)=\frac{f_j}{g_j},
\ees
where $f_j\in P$, $0\neq g_j\in P$, and the pairs $f_j,g_j$ are relatively prime for all $j$.  Notice that $e_j^d\in P_0$ and 
\bes
D(e_j^d)=S(e_j^d)=de_j^{d-1}\frac{f_j}{g_j}\in P_0. 
\ees
 Consequently, 
\bee\label{h9}
g_j|e_j^{d-1}
\eee
since the pair $f_j, g_j$ is relatively prime, that is, $g_j$ is a power of $e_j$.

If $d | \deg e_j$ then $e_j\in P_0$ and 
$S(e_j)=D(e_j)\in P_0$, i.e., we may assume that $g_j=1$. If $d\not|\deg e_j$ then there exist $1\leq i\leq n$ and $1\leq k<d$ such that $e_i=x_i\neq e_j$ and $x_i^ke_j\in P_0$. Then 
\bes
D(x_i^ke_j)=S(x_i^ke_j)=kx_i^{k-1}\frac{f_i}{g_i}e_j+x_i^k\frac{f_j}{g_j}\in P_0. 
\ees
Consequently, 
\bes
g_ig_j | kx_i^{k-1}f_ig_je_j+x_i^kf_jg_i. 
\ees
This implies that $g_j | x_i^kf_jg_i$. Since $f_j, g_j$ are relatively prime, we get $g_j | x_i^kg_i$.  By (\ref{h9}) $g_i$ is a power of $e_i=x_i$ and $g_j$ is a power of $e_j$. Since $i\neq j$, this   implies that $g_j\in K$ for all $j$. Consequently, $S(e_j)\in P$ and $S(P)\subseteq P$. 

Let us now show that the restriction of $S$ to $P$ is a Lie derivation,  i.e., 
\bee\label{h10}
S(\{u,v\})=\{S(u),v\}+\{u,S(v)\} 
\eee
 for all $u,v$ of the form (\ref{h5}). We prove (\ref{h10}) by induction on the polynomial length $d(u)+d(v)$. 
Suppose that $u=e_i$ and $v=e_j$. Since $e_i^d,e_j^d\in P_0$, we get 
\bes
S(\{e_i^d,e_j^d\})=D(\{e_i^d,e_j^d\})=\{D(e_i^d),e_j^d\}+\{e_i^d,D(e_j^d)\}\\
=\{S(e_i^d),e_j^d\}+\{e_i^d,S(e_j^d)\}=\{de_i^{d-1}S(e_i),e_j\}\nonumber de_j^{d-1}+de_i^{d-1}\{e_i,de_j^{d-1}S(e_j)\}\\
=d^2(d-1)e_i^{d-2}e_j^{d-1}S(e_i)\{e_i,e_j\}+d^2e_i^{d-1}e_j^{d-1}\{S(e_i),e_j\}\\+d^2(d-1)e_i^{d-1}e_j^{d-2}S(e_j)\{e_i,e_j\}+d^2e_i^{d-1}e_j^{d-1}\{e_i,S(e_j)\}. 
\ees
On the other hand, 
\bes
\{e_i^d,e_j^d\}=d^2e_i^{d-1}e_j^{d-1}\{e_i,e_j\}
\ees
and
\bes
S(\{e_i^d,e_j^d\})=S(d^2e_i^{d-1}e_j^{d-1}\{e_i,e_j\})\\
=d^2S(e_i^{d-1})e_j^{d-1}\{e_i,e_j\} 
+d^2e_i^{d-1}S(e_j^{d-1})\{e_i,e_j\}+d^2e_i^{d-1}e_j^{d-1}S(\{e_i,e_j\})\\
=d^2(d-1)e_i^{d-2}e_j^{d-1}S(e_i)\{e_i,e_j\}+d^2(d-1)e_i^{d-1}e_j^{d-2}S(e_j)\{e_i,e_j\}
+d^2e_i^{d-1}e_j^{d-1}S(\{e_i,e_j\}).
\ees
Comparing two values of $S(\{e_i^d,e_j^d\})$, we get 
\bes
S(\{e_i,e_j\})=\{S(e_i),e_j\}+\{e_i,S(e_j)\}.
\ees
Suppose that $d(v)\geq 2$ and $v=v_1v_2$. Then 
\bes
S(\{u,v\})=S(\{u,v_1v_2\})=S(v_1\{u,v_2\}+\{u,v_1\}v_2)\\
=S(v_1)\{u,v_2\}+v_1S(\{u,v_2\})+S(\{u,v_1\})v_2+\{u,v_1\}S(v_2). 
\ees
By the induction proposition, (\ref{h10}) is true for pairs $u,v_1$ and $u,v_2$, i.e., 
\bes
S(\{u,v_1\})=\{S(u),v_1\}+\{u,S(v_1)\}, \ \ \ S(\{u,v_2\})=\{S(u),v_2\}+\{u,S(v_2)\}. 
\ees
Then
\bes
S(\{u,v\})=S(v_1) \{u,v_2\}+v_1 \{S(u),v_2\}+v_1 \{u,S(v_2)\}\\
 +\{S(u),v_1\}v_2+\{u,S(v_1)\}v_2+\{u,v_1\}S(v_2)\\
= \{S(u),v_1v_2\}+\{u,S(v_1)v_2+v_1S(v_2)\}=\{S(u),v\}+\{u,S(v)\}. 
\ees
Consequently, $S$ is a derivation of a Poisson algebra $P$ and induces $D$ on $P_0$. $\Box$

\begin{lm}\label{m3} Every locally nilpotent derivation of the Poisson algebra $P_0$ is induced by a locally nilpotent $d$-derivation of the Poisson algebra $P=P\langle x_1,x_2,\ldots,x_n\rangle$. 
\end{lm}
\Proof 
Let $D$ be a locally nilpotent derivation of $P_0$ and let $S$ be a unique extension of $D$ to $P$. We have to show that $S$ is a locally nilpotent derivation of $P$. Notice that 
$$P_0\subset P$$
 is an integral extension of domains since $e_i^d\in P_0$ for all $i\geq 1$. 
According to a result of W.V. Vasconcelos \cite{Vasconcelos} (see also 
Proposition 1.3.37 from \cite[p. 41]{Es}), $S$ is locally nilpotent.
$\Box$\\

\section{Automorphisms of $P_0$}

\hspace*{\parindent}

As we noticed above, every $d$-graded automorphism of $P\langle x_1,x_2,\ldots,x_n\rangle$ induces an automorphism of $P_0$. In this section we prove the reverse of this statement for  $n>1$.

\begin{theor}\label{t1} Let $K$ be a field closed with respect to taking all $d$-roots of elements. Then every automorphism of  $P_0$ over $K$ is induced by a $d$-graded automorphism of $P=P\langle x_1,x_2,\ldots,x_n\rangle$ if $n>1$. 
\end{theor}
\Proof Let $\a$ be an automorphism of $P_0$. Denote the extension of $\a$ to the quotient field $Q(P_0)$ by the same symbol. We have $\frac{x_2}{x_1}\in Q(P_0)$. Suppose that 
\bee\label{h15}
\a\left(\frac{x_2}{x_1}\right)=\frac{f_2}{f_1}, 
\eee
where $f_1,f_2$ are relatively prime. 
Then 
\bes
\a\left(\frac{x_2^d}{x_1^d}\right)=\a\left(\frac{x_2}{x_1}\right)^d=\frac{f_2^d}{f_1^d}. 
\ees
Since $f_1,f_2$ are relatively prime it follows that $\a(x_1^d)=vf_1^d$ and $\a(x_2^d)=vf_2^d$ for some $v\in P$. 
Moreover, $\a(x_1^ix_2^{d-i})=vf_1^if_2^{d-i}$ for all $0\leq i < d$.

We have $vf_1^d, vf_2^d\in P_0$. If $K$ is a field of characteristic $p>0$ and $p$ divides $d$, then $f_1^d, f_2^d\in P_0$. Consequently, $v\in P_0$. 
Assume that $K$ is a field of charateristic $0$ or of characteristic $p>0$ and $p$ does not divide $d$. Let $\epsilon$ be a primitive $d$-root of unity. Consider the automorphism $\varepsilon$ of $Q(P)$ such that $\varepsilon(x_i)=\epsilon x_i$ for all $i$.  Notice that for any $f\in Q(P)$ we have $f\in Q(P_0)$ if and only if $\varepsilon(f)=f$. Then 
\bes
\frac{f_2}{f_1}=\varepsilon(\frac{f_2}{f_1})=\frac{\varepsilon(f_2)}{\varepsilon(f_1)} 
\ees
and $\varepsilon(f_1), \varepsilon(f_2)$ are relatively prime. Hence $f_1\varepsilon(f_2)=\varepsilon(f_1)f_2$ and, say, $f_1$ divides $\varepsilon(f_1)$ and $\varepsilon(f_1)$ divides $f_1$. 
Consequently $f_1$ and $\varepsilon(f_1)$ are proportional which is possible only if $f_1$ is a $d$-homogeneous element.  Similarly, $f_2$ is a $d$-homogeneous element. Then $f_1^d, f_2^d\in P_0$ and, consequently, $v\in P_0$.

This implies that 
\bes
x_1^d=\a^{-1}(v)\a^{-1}(f_1^d), \ \ x_2^d=\a^{-1}(v)\a^{-1}(f_2^d). 
\ees
Since $x_1^d$ is irreducible in $P_0$, this is possible only if $v\in K$. 
 Let $\mu\in K$ be a $d$-root of $v$, i.e., $\mu^d=v$. Replacing $f_1$ and $f_2$ by $\mu f_1$ and $\mu f_2$, we may assume that 
\bee\label{h16}
\a(x_1^d)=f_1^d, \ \ \a(x_2^d)=f_2^d. 
\eee
By (\ref{h15}) and (\ref{h16}), we get 
\bes
\a(x_1^{i_1}x_2^{i_2})=f_1^{i_1} f_2^{i_2},\,\textit{if}\,\,x_1^{i_1} x_2^{i_2}\in P_0.
\ees

Consider an arbitrary $e_i$ with $i\geq 3$. Suppose that $\deg e_i=s$. Then $y_i=\frac{e_i}{x_1^s}\in Q(P_0)$. Suppose that 
\bes
\a(y_i)=\a\left(\frac{e_i}{x_1^s}\right)=\frac{f_i}{g_i}, 
\ees
where $f_i,g_i$ are relatively prime. Then 
\bes
\a\left(\frac{e_i^d}{x_1^{sd}}\right)=\a\left(\frac{e_i}{x_1^s}\right)^d=\frac{f_i^d}{g_i^d}. 
\ees
Again $\a(e_i^d)=vf_i^d$ and $\a(x_1^{sd})=vg_i^d$ for some $v\in P$.  As above, we get that $f_i^d, g_i^d\in P_0, v\in K$, and we can assume that 
\bes
\a(e_i^d)=f_i^d, \  \ \a(x_1^{sd})=g_i^d. 
\ees
Then $f_1^{sd}=g_i^d$ and $g_i=\lambda f_1^s$, where $\lambda$ is a $d$-root of unity. After rescaling, we can assume that $g_i=f_1^s$ and
\bee\label{h17}
\a(e_i^d)=f_i^d, \  \ \a(y_i)=\a\left(\frac{e_i}{x_1^s}\right)=\frac{f_i}{f_1^s}, 
\eee
where $s=\deg e_i$ and $i\geq 3$. This is true for $i=2$ by (\ref{h15}) and (\ref{h16}).

Let $u=e_{i_1}\ldots e_{i_k}$ be an arbitrary element of $P$ of the form (\ref{h5}). We have 
\bee\label{h11}
u=x_1^s\frac{e_{i_1}}{x_1^{s_{i_1}}}\ldots \frac{e_{i_k}}{x_1^{s_{i_k}}}=x_1^s y_{i_1}\ldots y_{i_k},
\eee
where $s=s_{i_1}+\ldots+s_{i_k}$.  We have $d|s$ since $u\in P_0$. Then 
\bes
\a(u)=f_1^s\frac{f_{i_1}}{f_1^{s_1}}\ldots \frac{f_{i_k}}{f_1^{s_k}}=f_{i_1}\ldots f_{i_k}
\ees
by (\ref{h15}), (\ref{h16}), and (\ref{h17}).

Consequently, the polynomial endomorphism $\b$ of $P$, determined by $\b(e_i)=f_i$ for all $i\geq 1$, induces 
$\a$ on $P_0$. First we show that $\b$ is a polynomial automorphism of $P$. The elements (\ref{h4}) are algebraically independent and, consequently, the elements $e_1^d,\ldots,e_s^d,\ldots$ are algebraically independent. Since $\a$ is an automorphism and $\a(e_i^d)=f_i^d$ for all $i$ by (\ref{h16}) and (\ref{h17}), the elements $f_1^d,\ldots,f_s^d,\ldots$ are algebraically independent. Therefore the elements $f_1,\ldots,f_s,\ldots$ are algebraically independent and $\b$ is an injective endomorphism. Then $\b$ can be uniquely extended to an endomorphism of the quotient field  $P(x_1,x_2,\ldots,x_n)$ and we denote this extension also by $\b$. 

The restriction of $\b$ on $Q(P_0)$ is an automorphism since it coincides with the $\a$. Consider the space
\bes
V=Q(P_0)P\langle x_1,x_2,\ldots,x_n\rangle. 
\ees
By (\ref{h11}) every element $f\in P$ can be written as 
\bes
f=f_0+f_1x_1+\ldots +f_{d-1}x_1^{d-1}, 
\ees
where $f_0,f_1,\ldots,f_{d-1}\in K[t,y_2,\ldots,y_s,\ldots]$ and $t=x_1^d$. 
Hence $V$ is the $Q(P_0)$-span of the elements $1,x_1,x_1^2,\ldots,x_1^{d-1}$. 
If 
\bes
V=b_1Q(P_0)\oplus\ldots\oplus b_kQ(P_0), 
\ees
then 
\bes
\b(V)=\b(b_1)Q(P_0)+\ldots +\b(b_k)Q(P_0)
\ees
since $\b(Q(P_0))=Q(P_0)$. Notice that $\b(V)\subseteq V$. If $\b(V)\neq V$ then $\dim_{Q(P_0)}V<k$ and $\mathrm{Ker} \b \neq 0$. It is impossible for nonzero field endomorphisms. Consequently, $\b(V)= V$ and $e_i\in \b(V)$ for all $i$. Therefore $\b$ is an automorphism of the field $P(x_1,x_2,\ldots,x_n)$ and of the polynomial algebra $P\langle x_1,x_2,\ldots,x_n\rangle$. 

It remains to show that $\b$ is a Lie automorphism of $P$, i.e., 
\bee\label{h18}
\b(\{u,v\})=\{\b(u),\b(v)\}
\eee
for all $u,v$ of the form (\ref{h5}). We prove (\ref{h18}) by induction on the polynomial length $d(u)+d(v)$. 
Suppose that $u=e_i$ and $v=e_j$. Since $e_i^d,e_j^d\in P_0$, we get 
\bes
\b(\{e_i^d,e_j^d\})=\a(\{e_i^d,e_j^d\})=\{\a(e_i^d),\a(e_j^d)\}=\{\b(e_i^d),\b(e_j^d)\}\\
=\{\b(e_i)^d,\b(e_j)^d\}=d^2\b(e_i)^{d-1}\b(e_j)^{d-1}\{\b(e_i),\b(e_j)\}. 
\ees
On the other hand, 
\bes
\b(\{e_i^d,e_j^d\})=\b(d^2e_i^{d-1}e_j^{d-1}\{e_i,e_j\})=d^2\b(e_i)^{d-1}\b(e_j)^{d-1}\b(\{e_i,e_j\}). 
\ees
Comparing two values of $\b(\{e_i^d,e_j^d\})$, we get that (\ref{h18}) holds for $u=e_i$ and $v=e_j$.

Suppose that $d(v)\geq 2$ and $v=v_1v_2$. Then 
\bes
\b(\{u,v\})=\b(\{u,v_1v_2\})=\b(v_1\{u,v_2\}+\{u,v_1\}v_2) \\
=\b(v_1)\b(\{u,v_2\})+\b(\{u,v_1\})\b(v_2). 
\ees
By the induction proposition, we may assume that (\ref{h18}) is true for pairs $u,v_1$ and $u,v_2$. 
Then
\bes
\b(\{u,v\})=\b(v_1)\{\b(u),\b(v_2)\}+\{\b(u),\b(v_1)\}\b(v_2)\\
= \{\b(u),\b(v_1)\b(v_2)\}=\{\b(u),\b(v)\}. 
\ees
Consequently, $\b$ is an automorphism of $P$ and induces $\a$ on $P_0$. $\Box$

Let $\mathrm{Aut}_dP$ be the group of all $d$-graded automorphisms of the free Poisson algebra $P$.

\begin{co}\label{c1} Let $K$ be a field closed with respect to taking all $d$-roots of elements and let $E=\{\lambda\mathrm{id} | \lambda^d=1, \lambda\in K\}$, where $\mathrm{id}$ is the identity automorphism of $P$. Then 
\bes
\mathrm{Aut}\,P_0\cong  \mathrm{Aut}_dP/E.
\ees
\end{co}
\Proof Consider the homomorphism 
\bee
\psi: \mathrm{Aut}_dP\to{Aut}\,P_0
\eee
defined by $\psi(\a)=\overline{\a}$, where $\overline{\a}$ is the automorphism of $P_0$ induced by the $d$-graded automorphism $\a$ of $P$. 

By Theorem \ref{t1}, $\psi$ is an epimorphism. Let $\a \in \mathrm{Ker}\,\psi$. Then $\a(x_1)^d=x_1^d$. Consequently, 
 $\a(x_1)=\lambda x_1$ for some $d$th root of unity $\lambda\in K$. Extending $\a$ to $Q(P_0)$, we get $\a(x_i/x_1)=x_i/x_1$. 
Consequently, $\a(x_i)=\lambda x_i$ for all $i$ and $\a=\lambda \mathrm{id}$, i.e., $\a\in E$. Obviously, $E\subseteq \mathrm{Ker}\,\psi$. $\Box$

\section{Veronese subalgebras of polynomial algebras}

\hspace*{\parindent}

Let $A=K[x_1,x_2,\ldots,x_n]$ be the polynomial algebra over a field $K$ in the variables $x_1,x_2,\ldots,x_n$. Consider the 
grading 
\bes
A=A_0\oplus A_1\oplus \ldots \oplus A_{d-1}, 
\ees
where $d\geq 2$ and $A_i$ is the subspace of $A$ generated by all monomials of degree $kd+i$ for all $k\geq 0$.  This is a $\mathbb{Z}_d$-grading of $A$, i.e., $A_iA_j\subseteq A_{i+j}$ for all $i,j\in \mathbb{Z}_d$. The subalgebra $A_0$ is called the {\em Veronese subalgera of $A$ of degree $d$}.

\begin{co}\label{c2} Let $A=K[x_1,x_2,\ldots,x_n]$ be the polynomial algebra over a field $K$ of characteristic zero in $n\geq 2$ variables $x_1,x_2,\ldots,x_n$. Then every derivation of the Veronese subalgebra $A_0$ can be uniquely extended to a $d$-graded derivation of $K[x_1,x_2,\ldots,x_n]$. 
\end{co}

\begin{co}\label{c3} Let $A=K[x_1,x_2,\ldots,x_n]$ be the polynomial algebra over a field $K$ of characteristic zero in $n\geq 2$ variables $x_1,x_2,\ldots,x_n$. Then every locally nilpotent derivation of the Veronese subalgebra $A_0$ is induced by a locally nilpotent $d$-derivation of the polynomial algebra $K[x_1,x_2,\ldots,x_n]$. 
\end{co}

\begin{co}\label{c4} Let $A=K[x_1,x_2,\ldots,x_n]$ be the polynomial algebra in $n\geq 2$ variables $x_1,x_2,\ldots,x_n$ over a field $K$ closed with respect to taking all $d$-roots of elements. Then every automorphism of the Veronese subalgebra $A_0$ of degree $d$ is induced by a $d$-graded automorphism of $K[x_1,x_2,\ldots,x_n]$. 
\end{co}

This result is also proven in \cite{Kollar1}.

Let $\mathrm{Aut}_dA$ be the group of all $d$-graded automorphisms of the polynomial algebra $A$.

\begin{co}\label{c5} Let $K$ be a field closed with respect to taking all $d$-roots of elements and let $E=\{\lambda\mathrm{id} | \lambda^d=1, \lambda\in K\}$, where $\mathrm{id}$ is the identity automorphism of $A$. Then 
\bes
\mathrm{Aut}\,A_0\cong  \mathrm{Aut}_dA/E.
\ees
\end{co}

The proofs of Corollary \ref{c2}, Corollary \ref{c3}, Corollary \ref{c4}, and Corollary \ref{c5} repeat the polynomial parts of the proofs of Lemma \ref{m2}, Lemma \ref{m3}, Theorem \ref{t1}, and Corollary \ref{c1}, respectively.

Notice that these statements are not true for the polynomial algebra $A=K[x]$ in one variable $x$. In this case, the Veronese subalgebra $A_0$ of degree $d$ is the polynomial algebra in one variable $x^d$. Then the locally nilpotent derivation of $A_0$ determined by 
\bes
x^d\mapsto 1
\ees
cannot be induced by any derivation of $A$ and 
 the automorphism of $A_0$ determined by 
\bes
x^d\mapsto x^d+1
\ees
cannot be induced by any automorphism of $A$.

In addition, analogues of these results are not true for free associative algebras. In fact, if $B=K\langle x,y\rangle$ is the free associative algebra in the variables $x,y$ and $d=2$ then the Veronese subalgebra $B_0$ of degree $d$ is the free associative algebra in the variables $x^2,xy,yx,y^2$. It is easy to check that the locally nilpotent derivation of $B_0$ determined by 
\bes
x^2\mapsto 1, xy\mapsto 0, yx\mapsto 0, y^2\mapsto 0
\ees 
cannot be induced by any derivation of $B$ and the automorphism 
of $B_0$ determined by 
\bes
x^2\mapsto x^2+ 1, xy\mapsto xy, yx\mapsto yx, y^2\mapsto y^2 
\ees 
cannot be induced by any automorphism of $B$.

\bigskip

\begin{center}
	{\large Acknowledgments}
\end{center}

\hspace*{\parindent}

The second and third authors are grateful to Max-Planck Institute f\"ur Mathematik for
hospitality and excellent working conditions, where part of this work has been done.

The third author is supported by the grant of the Ministry of Education and Science of the Republic of Kazakhstan (project  AP14872073).

\end{document}